\documentclass[a4paper,10pt]{article}
\usepackage{amsmath,amsthm,amsfonts}
\usepackage{amssymb}
\usepackage{graphicx,epsfig}

\theoremstyle{definition}

\theoremstyle{remark}

\title{\bf Counting One-Vertex Maps\thanks{The authors acknowledge partial funding of this research by
ARRS of Slovenia, grants: L2-7230,P1-0294,L2-7207 and the PASCAL network of
excellence in the 6th framework. } }
\author{Alen Orbani\' c {\small\tt\ Alen.Orbanic@fmf.uni-lj.si},\\
Marko Petkov\v sek {\small\tt\ Marko.Petkovsek@fmf.uni-lj.si},\\
Toma\v{z} Pisanski\thanks{Joint position at the University of Primorska} {\small\tt\,\ Tomaz.Pisanski@fmf.uni-lj.si}, \ and\\
\medskip
Primo\v z Poto\v cnik {\small\tt\ Primoz.Potocnik@fmf.uni-lj.si}\\
Institute of Mathematics, Physics and Mechanics\\
University of Ljubljana\\
Jadranska 19, 1000 Ljubljana, SLOVENIJA}
\date{\today}
\begin{document}
\maketitle

\begin{abstract}
The number of distinct maps (pre-maps) with a single vertex and valence $d$ is computed for
any value of $d$.
The types of maps (pre-maps) that we consider depend on whether the underlaying graph (pre-graph) is 
signed or unsigned and directed or undirected.
\end{abstract}

\section{Introduction}
\label{section:introduction}

The motivation for this note lies in the fact that each orientable Cayley map of
valence $d$ is obtained via a regular covering construction from an orientable map
with one vertex, $k$ loops and $d-2k$ half-edges (see \cite{malnic}; for Cayley
maps in general, see for instance \cite{CayleyMaps}. For basic definitions of combinatorial maps 
see \cite{bookchapter}). In this
context, the problem of determining the number $\pi(d)$ of all non-isomorphic
one-vertex $d$-valent maps arises naturally.
For more detailed analysis of these embeddings see \cite{kwak}.
Coverings of graphs and pregraphs are combinatorially described in \cite{bicirculants}.

As we show in this paper, the number $\pi(d)$ equals the number of all
essentially distinct matchings in the complete graph $K_d$ with vertices
arranged as in the regular $d$-gon,
where two matchings are considered essentially the same whenever
one can be obtained from the other by a rotation or a reflection
of the $d$-gon. The latter can be obtained by the formula
\[
\pi(d)\ =\ \displaystyle \frac{1}{2d} \left(F(d) + R(d)\right)
\]
where
\[
F(d) = \left\{
\begin{array}{ll}
\displaystyle
\frac{d}{2}\left(f\left(\frac{d}{2}\right) + 2 f\left(\frac{d}{2} - 1\right)\right), & d {\rm\ even}, \\
\displaystyle
d\, f\left(\frac{d-1}{2}\right), & d {\rm\ odd},
\end{array}
\right.
\]
\[
f(n) \ =\ n! \sum_{0\le 2j\le n} \frac{2^{n-2j}}{(n-2j)!j!},
\]
\[
R(d)\ =\ \sum_{r|d} \varphi\left(\frac{d}{r}\right) \sum_{0\le 2j\le r} {r\choose 2j}(2j-1)!!
\left(\frac{d}{r}\right)^j w(r)^{r-2j},
\]
and
\[
w(r)\ = \left\{
\begin{array}{ll}
2, & \displaystyle r \,|\ (d/2) \\
1, & \displaystyle r \not|\,\ (d/2)
\end{array}
\right.
\]
(here we assume that $r \not|\,\ (d/2)$ when $d$ is odd).
Similar formulae are obtained for one-vertex maps without half-edges, non-orientable maps,
and directed maps.

We introduce a structure that we call a pre-graph; see also \cite{malnic}.
A {\em pre-graph} $G$ is a quadruple $G = (V,S,i,r)$ where $V$ is the set of
vertices, $S$ is the set of {\em arcs} (also known as semi-edges, darts, sides,
...), $i$ is the {\em initial mapping} $i:S \rightarrow V$, specifying the origin
or initial vertex for each arc, while $r$ is the {\em reversal involution}:
$r: S \rightarrow S$, $r^2 = 1$. We may also define the {\em terminal mapping} $t:S
\rightarrow V$ as $t(s) := i(r(s))$, specifying the terminal vertex for each arc.
An arc $s$
forms an edge $e = \{s,r(s)\}$, which is called {\em
proper} if $|e| = 2$ and is called a {\em half-edge} if $|e| = 1$. Define
$\partial(e) = \{i(s),t(s)\}$. A pre-graph without half-edges is called a {\em
(general) graph}. Note that $G$ is a graph if and only if the reversal involution has no
fixed points. A proper edge $e$ with $|\partial(e)| = 1$ is called a {\em loop}
and two edges $e, e'$ are {\em parallel} if $\partial(e) =
\partial(e')$. A graph without loops and parallel edges is called {\em simple}. The
{\em valence} of a vertex $v$ is defined as $val(v) = |\{s \in S| i(s) = v\}|$.
All pre-graphs in this note are connected unless stated otherwise.

Topologically, an oriented map is a 2-cell embedding of a graph into  
an orientable surface. However, in this paper we will operate with the 
following combinatorial description.
For us an oriented map is a triple
$(S,r,R)$ where $S$ is a finite non-empty set and $r$, $R$ are permutations on $S$
such that $r^2=1$ and $\langle r, R\rangle$ acts transitively on $S$ (see \cite{malnic}). Note that the vertices of the map correspond to the cycles in the cyclic decomposition of $R$. An isomorphism between two pre-maps $(S,r,R)$ and $(S',r',R')$ is any
bijection $\pi \colon S\to S'$ for which $\pi R = R'\pi$ and $\pi r = r'\pi$
holds. An automorphism of a map $(S,r,R)$ is thus a permutation of $S$ which
commutes with both $R$ and $r$. 


In addition to oriented maps we will also consider 
general (possibly non-orientable) (pre-)maps, which can be defined
as $M = (S,r,R,\lambda)$, where $S$, $r$ and $R$ are as in the definition of
oriented maps, and $\lambda$ is a {\em sign mapping} assigning 
either $1$ or $-1$ to each proper edge of the underlying graph of the map $M$.
Recall that each cycle $C = (s_1, \ldots, s_k)$ of $R$ corresponds to a vertex of the 
map. Substituting the cycle $C$ in $R$ with the reverse cycle and 
inverting  the $\lambda$-value of proper edges underlying the darts 
$s_1, \ldots, s_k$ other than loops, results in a new map $M'$,
which is said to be obtained from $M$ by a {\em local orientation change}.
Two general maps $M_1 = (S_1,r_1,R_1,\lambda_1)$ and $M_2 = (S_2,r_2,R_2, \lambda_2)$ are {\em isomorphic} if there exists a
general map $M' = (S_1, r_1, R', \lambda')$ obtained from $M_1$ by a series
of local orientation changes and a bijection 
$\pi: S_1 \to S_2$, such that
$\pi r_1 = r_2 \pi$, $\pi R_1 = R' \pi$ and $\pi \lambda_1 =  \lambda' \pi$.

\section{Counting one-vertex graphs and pre-graphs}

Let $p(d)$ denote the number of one-vertex pre-graphs of valence $d$. Since each
of them is determined by the number of loops, $p(d)$ can be computed using the
formula: $p(d) = 1 + \lfloor d/2 \rfloor$. This gives rise to the 
generating function $P(x) = 1/((1-x)^2(1+x))$.

If there are no pending edges, the situation becomes much simpler. Let $g(d)$
denote the number of one-vertex graphs. Then $g(d) = 0$, for $d$ odd, and $g(d) =
1$, for $d$ even. The corresponding generating function $G(x)$ is $G(x) =
1/(1-x^2)$.


\section{Counting one-vertex maps and pre-maps}

When counting (oriented) pre-maps the same pre-graph may give rise to more than one
pre-map.

Let $\pi(d)$  denote the number of oriented pre-maps whose underlying graphs 
are single
vertex pre-graphs, and let $\gamma(d)$  denote the number of oriented maps whose
underlying graphs are single vertex graphs.

Let $\pi_\tau(d)$  denote the number of non-isomorphic single vertex pre-maps of
type $\tau$ and valence $d$. The types of pre-maps that we consider are denoted by
$\bar S \bar D \bar G,
      S \bar D \bar G,
 \bar S      D \bar G,
      S      D \bar G,
 \bar S \bar D      G,
      S \bar D      G,
 \bar S      D      G$, and
 $     S      D      G$,
indicating whether the underlying pre-graphs
are signed or unsigned ($S$ resp.\ $\bar S$), directed or undirected ($D$ resp.\ $\bar D$),
graphs or pre-graphs ($G$ resp.\ $\bar G$). As it turns out, $\pi_\tau(d)$ and the various auxiliary
functions can be written in the same general form for all $\tau$, but with different
values of parameters (cf.\ Table \ref{stm}).
If necessary, we refer to the three symbols composing $\tau$ by $\tau_1$, $\tau_2$ and $\tau_3$.



Each oriented one-vertex pre-map is isomorphic to one of the form $(S,r,R)$ where
$S=\{1,\ldots,d\}$, $R=(1,2,\ldots,d)$. Such a pre-map can be represented by a
matching in the complete graph $K_d$ (possibly an empty one) in which two vertices
$i,j\in \{1,\ldots, d\}$ are matched whenever $r(i) = j$. Hence the number of all
one-vertex pre-maps $(S,r,R)$ with a given rotation $R$ is the same as the number
of all matchings (including the empty one) in $K_d$. This number is easily
computed to be
$$i(d) = \sum_{0 \leq 2k \leq d}{{d}\choose{2k}}(2k-1)!!,$$
where $(-1)!! = 1$.
Note that $i(d) = i(d-1) + (d-1)i(d-2)$ for all $d \ge 2$, and the exponential generating
function of this sequence is $\sum_{n=0}^\infty i(n) x^n/n! = \exp(x+x^2/2)$.

Of course, many of the above pre-maps are isomorphic. To compute the number of
non-isomorphic ones, let $I_d$ denote the set of all
permutations $r$ on $S$ with $r^2=1$, and
recall that two pre-maps $M=(S,r,R)$ and $M'=(S,r',R)$ are isomorphic if and only
if there exists a permutation on $S$ which centralizes $R$ and conjugates $r$ to $r'$.
Since the permutations of $S$ that centralize $R$ form the dihedral group $D_d$ of order $2d$,
it follows that the number of orientable non-isomorphic one-vertex pre-maps of valence $d$
equals the number of orbits of $D_d$ in its action on the set $I_d$ by conjugation.
If the pre-maps are represented by the matchings in $K_d$, as described above, then
the action of $D_d$ on $I_d$ corresponds to the natural action of $D_d$ on the set of all
matchings in $K_d$. In general, the number $\pi_\tau(d)$
of non-isomorphic one-vertex pre-maps of type $\tau$ and valence $d$ equals the number of orbits
of $D_d$ in its natural action on the set of matchings of type $\tau$ in $K_d$. The latter can be obtained
by the well-known Cauchy-Frobenius Lemma (also known as Burnside's Lemma):
\begin{equation}
\label{burnside}
\pi_\tau(d) \ =\ \frac{1}{|D_d|} \sum_{\sigma \in D_d} |{\rm Fix}_\tau (\sigma)|
\end{equation}
where ${\rm Fix}_\tau (\sigma)$ denotes the set of matchings of type $\tau$ in $K_d$ invariant under $\sigma$.
For more information on this method, see \cite{burnsideSchat}.

\subsection{Fixed points of reflections}

In order to compute the sum in (\ref{burnside}), assume first that $\sigma$ is a reflection. We distinguish two cases.

If $d$ is even there are two types of reflections: either across a median or across a main diagonal.
Let $\sigma$ be the reflection across a median, and let $L$ denote the set of $n = d/2$ vertices of $K_d$
on one side of the median. For each $u \in L$, denote by $u'$ its mirror image across the median. We will derive a recurrence
satisfied by $f_\tau(n) := |{\rm Fix}_\tau (\sigma)|$, using the so-called {\em method of distinguished element}.
Assume that $n \ge 2$,
pick any vertex $u \in L$, and partition ${\rm Fix}_\tau (\sigma)$ into subsets $A$ and $B$ where $A$ contains
those matchings in which $u$ is matched with $u'$ or is left unmatched, and $B$ contains those matchings in which
$u$ is matched with $v$ or $v'$ where $v$ is one of the remaining $n-1$ vertices in $L$. Denote by $s_\tau$ the number of ways
in which $u$ can be matched with $u'$ (including leaving it unmatched), and by $t_\tau$ the number of ways in which $u$ can be
matched with $v$. Because of symmetry, the number of ways in which $u$ can be matched with $v'$ is also $t_\tau$. The values
of $s_\tau$ and $t_\tau$ depend on the type $\tau$ of the problem considered, and are shown in Table~\ref{stm}.
\begin{table}[h]
\[
\begin{array}{l|cccccccc}
\tau &
\bar S \bar D \bar G &
     S \bar D \bar G &
\bar S      D \bar G &
     S      D \bar G &
\bar S \bar D      G &
     S \bar D      G &
\bar S      D      G &
     S      D      G \\
\hline\\*[-9pt]
s_\tau &     2 & 3 & 1 & 1 & 1 & 2 & 0 & 0 \\
t_\tau &     1 & 2 & 2 & 4 & 1 & 2 & 2 & 4 \\
m_\tau &     2 & 3 & 3 & 5 & 1 & 2 & 2 & 4 \\
\end{array}
\]
\[
\begin{array}{lll}
s_\tau & \ldots & {\rm \ the\ number\ of\ ways\ to\ match\ } u \in L {\rm \ with\ } u' \\
t_\tau & \ldots & {\rm \ the\ number\ of\ ways\ to\ match\ } u \in L {\rm \ with\ } v \in L \setminus\{u\} \\
m_\tau & \ldots & {\rm \ the\ number\ of\ ways\ to\ match\ the\ two\ vertices\ on\ the\ mirror}
\end{array}
\]
\caption{The values of parameters $s_\tau, t_\tau, m_\tau$ for the types of pre-maps considered}
\label{stm}
\end{table}
Then $|A| = s_\tau f_\tau(n-1)$ and $|B| = 2 t_\tau (n-1)\, f_\tau(n-2)$, hence
\begin{equation}
\label{rec}
f_\tau(n) = s_\tau\, f_\tau(n-1) + 2 t_\tau\, (n-1) f_\tau(n-2), \quad {\rm for\ } n \ge 2,
\end{equation}
with $f_\tau(0)=1$, $f_\tau(1)=s_\tau$.
To solve (\ref{rec}), let $G_\tau(x) = \sum_{n=0}^\infty f_\tau(n) x^n/n!$ be the exponential
generating function of the sequence $\langle f_\tau(n)\rangle_{n=0}^\infty$. Then it follows from (\ref{rec})
and the initial values that $G_\tau(x)$ satisfies the differential equation
\[
G_\tau'(x) = (s_\tau + 2\, t_\tau\,x) G_\tau(x), \quad G_\tau(0) = 1,
\]
whence
\begin{equation}
\label{gf}
\sum_{n=0}^\infty f_\tau(n) \frac{x^n}{n!}\ =\ \exp(s_\tau\,x + t_\tau\,x^2).
\end{equation}
By expanding the right-hand side into power series and comparing coefficients we find the solution
\begin{equation}
\label{expl}
f_\tau(n) \ =\ n! \sum_{0\le 2j\le n} \frac{s_\tau^{n-2j} t_\tau^j}{(n-2j)!j!}
     \ =\ \sum_{0\le 2j\le n} s_\tau^{n-2j} \,(2 t_\tau)^j {n\choose 2j}(2j-1)!!
\end{equation}
where $0^0 = 1$. Note that (\ref{expl}) can also be obtained by a counting argument:
To construct a matching $M$ which is invariant under $\sigma$, select $2j$ vertices from among the $n$ vertices in $L$,
then construct a perfect matching on these $2j$ vertices. This can be done in ${n\choose 2j}(2j-1)!!$ ways. As above,
there are $2t_\tau$ ways to match the two elements in each of the $j$ pairs, yielding the factor $(2t_\tau)^j$, and $s_\tau$ ways
to match each of the remaining $n-2j$ vertices to its mirror image, yielding the factor $s_\tau^{n-2j}$.

By comparing (\ref{gf}) to the generating function of Hermite polynomials
\[
\sum_{n=0}^\infty H_n(z) \frac{x^n}{n!}\ =\ \exp(2z\,x - x^2)
\]
we can also express $f_\tau(n)$ in terms of the $n$-th Hermite polynomial as
\begin{equation}
\label{herm}
f_\tau(n) = \left(i \sqrt{t_\tau}\right)^n H_n\left(\frac{s_\tau}{2i \sqrt{t_\tau}}\right).
\end{equation}

If $\sigma$ is the reflection across a main diagonal then $|{\rm Fix}_\tau (\sigma)|\ =\ m_\tau\, f(d/2 - 1)$
where $m_\tau$ is the number of ways in which it is possible to match the two vertices on the mirror
with each other. The value of $m_\tau$ depends on the type $\tau$
of the problem considered, and is shown in Table \ref{stm}.

If $d$ is odd there is only one type of reflections, and $|{\rm Fix}_\tau (\sigma)|\ =\ f_\tau((d-1)/2)$
for pre-maps and $0$ for maps.
Thus the total contribution $F_\tau(d)$ of the $d$ reflections to the sum in (\ref{burnside}) is
\begin{equation}
\label{totalRefl}
F_\tau(d) = \left\{
\begin{array}{ll}
\displaystyle
\frac{d}{2}\left(f_\tau\left(\frac{d}{2}\right) + m_\tau\, f_\tau\left(\frac{d}{2} - 1\right)\right), & d {\rm\ even}, \\
\displaystyle
d\, f_\tau\left(\frac{d-1}{2}\right), & d {\rm\ odd\ and\ }\tau_3 = \bar G, \\
0, & d {\rm\ odd\ and\ }\tau_3 = G,
\end{array}
\right.
\end{equation}
where $f_\tau$ is given by any of (\ref{rec}), (\ref{gf}), (\ref{expl}), or (\ref{herm}).

\subsection{Fixed points of rotations}

Now assume that $\sigma$ is
the counter-clockwise rotation of $2\pi k_\sigma / d$ where $0 \le k_\sigma < d$.
In how many ways can we construct a matching $M$ of $K_d$ which is invariant under $\sigma$?

Let $r = \gcd(d, k_\sigma)$. Then $\sigma$ has $r$ orbits in $V(K_d)$, each containing
$d/r$ vertices. Let $C$ denote a set of $r$ consecutive vertices of $K_d$. Since $C$
contains one representative from each orbit, it suffices to define $M$ on $C$, and
to extend it to $V(K_d)\setminus C$ by symmetry. Hence we can also
think of $M$ as a matching of orbits. Assume that $2j$ of the $r$ orbits are matched in pairs,
while the rest remain unmatched or are matched with themselves (the latter is possible only
if antipodal vertices belong to the same orbit, i.e., if $d$ is even and $r \,|\, d/2$).
There are ${r\choose 2j}$ ways to select the $2j$ orbits,
and $(2j-1)!!$ ways to group them into pairs. In each of the $j$ pairs of orbits $(\alpha_i,
\beta_i)$, $i = 1, 2, \ldots, j$, the vertex in $\alpha_i \cap C$ can be matched with any
of the $d/r$ vertices in $\beta_i$ in $t_\tau$ ways, and each of the remaining $r-2j$
orbits can be matched to themselves (or be left unmatched) in $w_\tau(r)$ ways where
\[
w_\tau(r)\ = \left\{
\begin{array}{ll}
s_\tau, & \displaystyle r \,|\ (d/2) \\
0, & \displaystyle r \not|\,\ (d/2) {\rm\ \;and\ } \tau_3 = G \\
1, & \displaystyle r \not|\,\ (d/2) {\rm\ \;and\ } \tau_3 = \bar G
\end{array}
\right.
\]
(for the values of $s_\tau$ and $t_\tau$, see Table \ref{stm}).
Now for each divisor $r$ of $d$, there are $\varphi(d/r)$ rotations $\sigma$ in $D_d$
having $\gcd(d, k_\sigma) = r$. Hence the total contribution $R_\tau(d)$ of the $d$ rotations to the sum in
(\ref{burnside}) is
\begin{equation}
\label{totalRot}
R_\tau(d)\ =\ \sum_{r|d} \varphi\left(\frac{d}{r}\right) \sum_{0\le 2j\le r} {r\choose 2j}(2j-1)!!
\left(\frac{t_\tau d}{r}\right)^j w_\tau(r)^{r-2j}
\end{equation}
where, as before, $0^0 = 1$.

\subsection{The master formula}

From (\ref{burnside}) it follows that the number of non-isomorphic single vertex pre-maps of
valence $d$ is
\begin{equation}
\label{master}
\pi_\tau(d)\ =\ \displaystyle \frac{1}{2d} \left(F_\tau(d) + R_\tau(d)\right)
\end{equation}
where $F_\tau(d)$ resp.\ $R_\tau(d)$ are given by (\ref{totalRefl}) resp.\ (\ref{totalRot}),
and the values of parameters $s_\tau, t_\tau, m_\tau$ for each type $\tau$ of pre-maps considered are given
in Table \ref{stm}.

\section{Additional formul\ae\ and tables}

Some of the sequences encountered in this paper can be found in the
{\em The Online Encyclopedia of Integer Sequences\/} (OEIS, \cite{OEIS}).

\begin{table}[h]
\begin{center}
\begin{tabular}{c|c|c}
sequence & OEIS ID number & exponential generating function  \\
\hline\\*[-9pt]
$\langle f_{\bar S \bar D \bar G}(n)\rangle_{n=0}^\infty$    &  A000898   & $\exp(x^2+2x)$   \\
$\langle f_{     S \bar D \bar G}(n)\rangle_{n=0}^\infty$    &            & $\exp(2x^2+3x)$   \\
$\langle f_{\bar S      D \bar G}(n)\rangle_{n=0}^\infty$    &  A115329   & $\exp(2x^2+ x)$   \\
$\langle f_{     S      D \bar G}(n)\rangle_{n=0}^\infty$    &            & $\exp(4x^2+ x)$   \\
\hline\\*[-9pt]
$\langle f_{\bar S \bar D      G}(n)\rangle_{n=0}^\infty$    &  A047974   & $\exp(x^2+ x)$   \\
$\langle f_{     S \bar D      G}(n)\rangle_{n=0}^\infty$    &            & $\exp(2x^2+2x)$   \\
$\langle f_{\bar S      D      G}(n)\rangle_{n=0}^\infty$    &            & $\exp(2x^2)$   \\
$\langle f_{     S      D      G}(n)\rangle_{n=0}^\infty$    &            & $\exp(4x^2)$   \\
\hline\\*[-9pt]
$\langle f_{\bar S      D      G}(2n-2)\rangle_{n=1}^\infty$ &  A052714   & $(1 - \sqrt{1 - 8x})/4$   \\
$\langle f_{     S      D      G}(2n-2)\rangle_{n=1}^\infty$ &  A052734   & $(1 - \sqrt{1 - 16x})/8$   \\
\hline\\*[-9pt]
$\langle \pi_{\bar S \bar D    G}(2n)\rangle_{n=0}^\infty$   &  A054499   &                  \\

\end{tabular}
\caption{ID numbers and generating functions of some of our sequences}
\end{center}
\end{table}

When $s_\tau = 0$, the formula giving $f_\tau(d)$ can be expressed in closed form. Thus, for $d$ even,
\begin{eqnarray*}
f_{\bar S D G}(d) &=& 2^d (d - 1)!!, \\
f_{     S D G}(d) &=& (2 \sqrt 2)^d (d - 1)!!.
\end{eqnarray*}
Also, when $w_\tau(r) = 0$ for all $r$, the double sum in the formula giving $R_\tau(d)$ reduces to a single sum.
For $d$ even we thus have
\begin{eqnarray*}
R_{\bar S D G}(d) &=& \sum_{r|d,\ r{\rm\,even}} \varphi\left(\frac{d}{r}\right) (r-1)!!
\left(\frac{2d}{r}\right)^{r/2}, \\
R_{     S D G}(d) &=& \sum_{r|d,\ r{\rm\,even}} \varphi\left(\frac{d}{r}\right) (r-1)!!
\left(\frac{4d}{r}\right)^{r/2}.
\end{eqnarray*}

In Tables \ref{premaptable} resp.\ \ref{maptable} we list the numbers of
non-isomorphic single vertex pre-maps resp.\ maps of valence $d$ for small values
of $d$. For instance, in \cite{watkins} the six pre-maps of valence five are
discussed in detail.

\begin{table}[h]
\begin{center}
\begin{tabular}{r|rrrr}
$d$ & $\pi_{\bar S \bar D \bar G}(d)$ & $\pi_{S \bar D \bar G}(d)$ & $\pi_{\bar S D \bar G}(d)$ & $\pi_{S D \bar G}(d)$ \\
\hline\\*[-9pt]
1  &  1         & 1            & 1             & 1               \\
2  &  2         & 3            & 2             & 3               \\
3  &  2         & 3            & 2             & 3               \\
4  &  5         & 11           & 6             & 14              \\
5  &  6         & 15           & 11            & 33              \\
6  &  17        & 60           & 37            & 167             \\
7  &  27        & 125          & 100           & 619             \\
8  &  83        & 529          & 405           & 3686            \\
9  &  185       & 1663         & 1527          & 18389           \\
10 &  608       & 7557         & 6824          & 120075          \\
11 &  1779      & 31447        & 30566         & 706851          \\
12 &  6407      & 155758       & 151137        & 5032026         \\
13 &  22558     & 763211       & 757567        & 33334033        \\
14 &  87929     & 4089438      & 4058219       & 255064335       \\
15 &  348254    & 22190781     & 22150964      & 1855614411      \\
16 &  1456341   & 127435846    & 127215233     & 15129137658     \\
17 &  6245592   & 745343353    & 745057385     & 119025187809    \\
18 &  27766356  & 4549465739   & 4547820514    & 1026870988199   \\
19 &  126655587 & 28308456491  & 28306267210   & 8640532108675   \\
20 &  594304478 & 182435301597 & 182422562168  & 78446356190934  \\
\end{tabular}
\caption{The numbers of non-isomorphic one-vertex pre-maps }
\label{premaptable}
\end{center}
\end{table}

\begin{table}[h]
\begin{center}
\begin{tabular}{r|rrrr}
$d$ & $\pi_{\bar S \bar D G}(d)$ & $\pi_{S \bar D G}(d)$ & $\pi_{\bar S D G}(d)$ & $\pi_{S D G}(d)$ \\
\hline\\*[-9pt]
2  &  1                & 2                   & 1                    & 2                     \\
4  &  2                & 6                   & 3                    & 9                     \\
6  &  5                & 26                  & 13                   & 90                    \\
8  &  17               & 173                 & 121                  & 1742                  \\
10 &  79               & 1844                & 1538                 & 48580                 \\
12 &  554              & 29570               & 28010                & 1776358               \\
14 &  5283             & 628680              & 618243               & 79080966              \\
16 &  65346            & 16286084            & 16223774             & 4151468212            \\
18 &  966156           & 490560202           & 490103223            & 250926306726          \\
20 &  16411700         & 16764409276         & 16761330464          & 17163338379388        \\
22 &  312700297        & 639992710196        & 639968394245         & 1310654311464970      \\
24 &  6589356711       & 26985505589784      & 26985325092730       & 110531845060209836    \\
\end{tabular}
\caption{The numbers of non-isomorphic one-vertex maps}
\label{maptable}
\end{center}
\end{table}

Using methods of \cite{disconstruct} one can easily extend the counting 
to graphs with one-vertex connected components. Motivated by \cite{kwak}, it would be worthwhile to
extend this analysis to dipoles or any two-vertex graphs or pre-graphs.



\end{document}